\documentclass[a4paper,11pt]{article}
\usepackage[top=2.5cm,bottom=2.5cm,left=2.2cm,right=2.2cm]{geometry}
\usepackage{amsfonts}
\usepackage{mathrsfs,amscd,amssymb,amsthm,amsmath,bm,graphicx,psfrag,subfigure,url,mathtools}
\usepackage{pict2e}
\usepackage{stfloats}
\usepackage{psfrag,amsmath}
\usepackage{tikz}
\usepackage{indentfirst}
\usepackage{hyperref}
\usepackage{bookmark}
\usepackage{enumerate}
\usepackage{latexsym,euscript,epic,eepic,color}
\usepackage{multirow}
\usepackage{multicol}
\usepackage{longtable}
\usepackage{adjustbox}
\usepackage[all]{xy}
\usepackage{setspace}
\usepackage{epstopdf}
\allowdisplaybreaks
\usepackage{authblk}
\usepackage{cmap}
\usepackage{pifont}

\makeatletter

\renewcommand{\@seccntformat}[1]{{\csname the#1\endcsname}{\normalsize .}\hspace{.5em}}
\makeatother

\usepackage{ifpdf}
\usepackage{indentfirst}

\def \[{\begin{equation}}
\def \]{\end{equation}}

\usepackage{epstopdf}
\numberwithin{equation}{section}
\newtheorem{theorem}{Theorem}[section]
\newtheorem{corollary}[theorem]{Corollary}
\newtheorem{definition}[theorem]{Definition}
\newtheorem{conjecture}[theorem]{Conjecture}

\newtheorem{problem}{Problem}

\allowdisplaybreaks

\newenvironment{wst}
{\setlength{\leftmargini}{1.5\parindent}
 \begin{itemize}
 \setlength{\itemsep}{-1.1mm}}
{\end{itemize}}

\begin{document}
\baselineskip=0.23in

\title{A survey of edge-spectral-Tur\'an type problems in spectral graph theory: Results, conjectures and open problems}%\thanks{%X.G. is financially supported by the National Natural Science Foundation of China (Grant Nos. 12171190, 11671164)\\[3pt] 
%\hspace*{5mm}{\it Email addresses}: ytyumath@sina.com (Y. Yu)}}

%\author[1]{Xianya Geng}% Author
%\author[2]{Jing Huang}
\author[,1]{Yuantian Yu\thanks{ytyumath@sina.com (Y. Yu)}}
\author[,2]{Huihui Zhang\thanks{zhanghhmath@163.com (H. Zhang)}}
\author[,3]{Minjie Zhang\thanks{Corresponding author, zhangmj1982@qq.com (M. Zhang)}}
%Communication author
%\affil[1]{School of Mathematics and Big date, Anhui University of Science and Technology, Huainan 232001, China}
%\affil[2]{School of Mathematics and Information Science, Guangzhou University, Guangzhou, 510006, China}% bit
\affil[1]{School of Science, East China University of Technology, Nanchang 330013, China}
\affil[2]{School of Mathematics, Luoyang Normal University, Luoyang 471000, China}
\affil[3]{School of Mathematics and Statistics, \linebreak Hubei University of Arts and Science, Xiangyang 441053, China}
\date{\today}
\maketitle

\begin{abstract}
The edge-spectral-Tur\'an type problem is also called the Brualdi-Hoffman-Tur\'an type problem, which is a central topic in spectral graph theory, seeking to determine the maximum spectral radius $\lambda(G)$ of an $F$-free graph $G$ with $m$ edges. This problem has attracted significant attention in recent years. In this paper, we will sort out several closely related results in this type of problem and then propose some conjectures for further research.

\vspace{2mm}

\end{abstract}

\section{\normalsize Introduction}
In this survey, we consider graphs that are finite, undirected, and simple (i.e., no loops or multiple edges). Unless explicitly stated otherwise, we adhere to the standard notation and terminology commonly found in graph theory (see, for example, Bollob\'as \cite{Bollo1}, Godsil and Royle \cite{Godsil1}). 

Let $G =(V(G), E(G))$ be a graph with vertex set $V(G) = \{v_1,v_2, \dots, v_n\}$ and edge set $E(G)$. Two vertices $v_i$ and $v_j$ in $G$ are said to be \textit{adjacent} (or \textit{neighbours}) if they are joined by an edge, denoted by \( v_i \sim v_j \). The adjacency matrix $A(G)=(a_{ij})$ of $G$ is an $n\times n$\ $(0,1)$-matrix with $a_{ij}=1$ if and only if $v_i,v_j$ are adjacent. As $A(G)$ is real and symmetric, all the eigenvalues $\lambda_1(G)\geqslant\lambda_2(G)\geqslant\cdots\geqslant\lambda_n(G)$ of $A(G)$ are real. As usual, the \textit{spectral radius} $\lambda(G)$ of graph $G$ is the largest eigenvalue of $A(G)$ in modulus. By the Perron-Frobenius theorem, $\lambda(G)=\lambda_1(G)$ and there exists a positive eigenvector $\boldsymbol{x}$ corresponding to $\lambda(G)$, known as the Perron vector of $G$, if $G$ is connected.

Given a family of graphs $ \mathcal{F} $, a graph $ G $ is said to be $\mathcal{F} $-free if it contains no subgraph isomorphic to any member of $ \mathcal{F} $. If $ \mathcal{F} $ consists of a single graph $ F $, $ G $ is referred to as $ F $-free. Let $ \mathcal{G}(m,  \mathcal{F}) $ denote the set of $  \mathcal{F} $-free graphs with $ m $ edges and no isolated vertices. In spectral graph theory, there are two types of extremal problems. One is the vertex-spectral-Tur\'an type problem, also called the Brualdi-Solheid-Tur\'an type problem \cite{BS1986}:
\begin{problem}[The vertex-spectral-Tur\'an type problem]\label{problem01}
Determine the maximum spectral radius among all $ \mathcal{F} $-free graphs on $n$ vertices and identify the corresponding extremal graphs.
\end{problem}
The other is the edge-spectral-Tur\'an type problem, also called the Brualdi-Hoffman-Tur\'an type problem \cite{Brualdi-Hoffman}:
\begin{problem}[The edge-spectral-Tur\'an type problem]\label{problem1}
Determine the maximum spectral radius among all graphs in $\mathcal{G}(m, \mathcal{F})$ and identify the corresponding extremal graphs.
\end{problem}

Li, Liu and Feng \cite{Li-Liu-Feng} surveyed a large number of vertex-spectral-Tur\'an type results. In past five years, more and more nice edge-spectral-Tur\'an type results have emerged. In this paper, we survey them in detail. In order to do so, we introduce some necessary notations and terminologies first.

Let $G$ be a graph, and let $u,v$ be two vertices of $G$. If $u\sim v$, let $G - uv$ (resp. $G - u$) denote the graph obtained from $G$ by deleting the edge $uv$ (resp. the vertex $u$); this notation extends naturally to the deletion of multiple edges or vertices. Similarly, if $u\nsim v$, let $G + uv$ denote the graph obtained from $G$ by adding an edge between $u$ and $v$. The set of neighbors of $v$ in $G$ is denoted by $N_G(v)$; its size, called the \textit{degree} of $v$ in $G$, is denoted by $d_G(v)$. Let $N_G[v] = N_G(v) \cup \{v\}$ (the closed neighborhood of $v$). 

As usual, let $P_n,\,C_n$ and $K_n$ be the path, the cycle and the complete graph of order $n$, respectively. And let $K_{a,b}$ be the complete bipartite graph with parts of sizes $a$ and $b$, respectively. Let $G_1$ and $G_2$ be two graphs, define $G_1\cup G_2$ to be their disjoint union. Then $G_1\vee G_2$ is defined to be their \textit{join}, obtained from $G_1\cup G_2$ by adding edges to connect each vertex of $G_1$ with all vertices of $G_2$. 

A graph is said to be \textit{properly coloured} if each vertex is assigned a color such that the two endpoints of every edge have distinct colors. The \textit{chromatic number} \(\chi(G)\) of a graph $G$ is defined as the minimum number $k$ for which $G$ can be properly colored using $k$ colors. A graph $G$ is called \textit{color-critical} if there exists an edge $e$ in $G$ such that the chromatic number of $G - e$ satisfies $\chi(G - e) < \chi(G)$.

\section{\normalsize Forbidding subgraphs involving consecutive odd cycles}
The study of Problem \ref{problem1} dates back to 1970, when Nosal~\cite{Nosal} obtained the following result.
\begin{theorem}[Nosal~\cite{Nosal}, 1970]\label{thm:2.1}
For all graphs $G \in \mathcal{G}(m, C_3),$ one has $\lambda(G) \leqslant \sqrt{m}$.
\end{theorem} 

Nikiforov~\cite{VNK2009} identified the corresponding extremal graphs in Theorem~\ref{thm:2.1}.

\begin{theorem}[Nikiforov~\cite{VNK2009}, 2009]\label{thm:2.2}
For all graphs $G \in \mathcal{G}(m, C_3),$ one has $\lambda(G) \leqslant \sqrt{m}$. Equality holds if and only if $G$ is a complete bipartite graph.
\end{theorem} 
Let \( G \) be a graph. A \textit{blow-up} of \( G \) is a new graph obtained from \( G \) by replacing each vertex \( x \in V(G) \) with an independent set \( I_x \), in which, for any two vertices \( x, y \in V(G) \), we add all edges between \( I_x \) and \( I_y \) if \( x\sim y \). Lin, Ning and Wu \cite{9} extended Theorem~\ref{thm:2.2} and obtained the following two results.

\begin{theorem}[Lin, Ning and Wu \cite{9}, 2021]\label{thm:5.2}
For all graphs $G \in \mathcal{G}(m, C_3)$ with $m\geqslant 2$, one has  
$$
\lambda_1^2(G) + \lambda_2^2(G) \leqslant m.
$$
Equality holds if and only if \( G \) is a blow-up of some member of \( \mathcal{G} \) in which  
$$
\mathcal{G} = \{ P_2,\ 2P_2,\ P_4,\ P_5 \}.
$$ 
\end{theorem}
\begin{theorem}[Lin, Ning and Wu~\cite{9}, 2021]\label{thm:2.3}
Let $G$ be a non-bipartite graph in $\mathcal{G}(m, C_3)$, one has $\lambda(G) \leqslant \sqrt{m - 1}$. Equality holds if and only if $G \cong C_5$.
\end{theorem}

Note that the bound in Theorem~\ref{thm:2.3} is attainable only when $m=5.$ For convenience, let $S_k(K_{a,b})$ denote the graph obtained from the complete bipartite graph $K_{a,b}$ by replacing one of its edges by a path of length $k+1,$ where $k\geq1$ and $b\geqslant a\geqslant1.$ Zhai and Shu~\cite{26} obtained the following result, which extended both Theorem~\ref{thm:2.2} and Theorem~\ref{thm:2.3}.
\begin{theorem}[Zhai and Shu~\cite{26}, 2022]\label{thm:2.4}
Let $G$ be a non-bipartite graph in $\mathcal{G}(m, C_3)$, one has $\lambda(G) \leqslant \lambda(S_1(K_{2, \frac{m - 1}{2}}))$. Equality holds if and only if $m$ is odd and $G \cong S_1(K_{2, \frac{m - 1}{2}}).$
\end{theorem}

The bound in Theorem~\ref{thm:2.4} is attainable only when $m$ is odd. Zhai and Shu~\cite{26} posed the following problem:
\begin{problem}[Zhai and Shu~\cite{26}, 2022]\label{pb001}
For even $m,$ what is the extremal graph attaining the maximum spectral radius over all $C_3$-free non-bipartite graphs with $m$ edges?
\end{problem}
%Recall that, for \( k \geq 1 \) and \( b \geq a \geq 1 \), \( S_k(K_{a,b}) \) is the graph obtained from \( K_{a,b} \) by replacing one of its edges by a path of length \( k + 1 \). 
For odd \( k \) and even \( m \), let \( S_k^+(K_{2, \frac{m-k-1}{2}}) \) be the graph obtained from \( S_k(K_{2, \frac{m-k-1}{2}}) \) by attaching a pendant edge to a vertex of maximum degree; if \( m - k \equiv 2 \pmod{3} \), then let \( S_k^-(K_{3, \frac{m-k+1}{3}}) \) be the graph obtained from \( K_{3, \frac{m-k+1}{3}} \) by the following two steps: First delete one edge of \( K_{3, \frac{m-k+1}{3}} \), which gives a \( 2 \)-degree vertex, say \( v \), in the resulting graph, say \( K_{3, \frac{m-k+1}{3}}^- \). Then replacing an edge incident with \( v \) of \( K_{3, \frac{m-k+1}{3}}^- \) by a path of length \( k + 1 \). Li, Feng and Peng \cite{Li-Feng-Peng} solved Problem~\ref{pb001} for even $m\geqslant 4.7\times 10^5$. Recently, Li, Zhao and Shen \cite{LZS} improved the result in \cite{Li-Feng-Peng}, solved Problem~\ref{pb001} for even $m\geqslant 70$.
\begin{theorem}[Li, Zhao and Shen~\cite{LZS}, 2025+]\label{thm:2.7}
Let $m\geqslant 70$ be even, and let $G$ be a non-bipartite graph in $\mathcal{G}(m, C_3)$.
\begin{wst}
\item[{\rm (a)}] If $m\equiv0 \pmod{3}$, then $\lambda(G)\leq\lambda(S_1^-(K_{3, \frac{m}{3}}))$ with equality if and only if $G\cong S_1^-(K_{3, \frac{m}{3}})$.
\item[{\rm (b)}] If $m\equiv1 \pmod{3}$, then $\lambda(G)\leq\lambda(S_1(K_{3, \frac{m-1}{3}}))$ with equality if and only if $G\cong S_1(K_{3, \frac{m-1}{3}})$.
\item[{\rm (c)}] If $m\equiv2 \pmod{3}$, then $\lambda(G)\leq\lambda(S_1^+(K_{2, \frac{m-2}{2}}))$ with equality if and only if $G\cong S_1^+(K_{2, \frac{m-2}{2}})$.
\end{wst}
\end{theorem}

Observe that $C_5$ is the unique extremal graph in Theorem~\ref{thm:2.3}, and all of the extremal graphs in Theorems~\ref{thm:2.4} and \ref{thm:2.7} contain $C_5$ as a subgraph. Sun and Li \cite{sl2023} extended Theorem \ref{thm:2.3} as follows.
\begin{theorem}[Sun and Li~\cite{sl2023}, 2023]\label{thm:2.5}
Let $G$ be a non-bipartite graph in $\mathcal{G}(m, \{C_3, C_5\})$, one has
\[\label{e1}
\lambda^4(G) \leqslant \sum_{v \in V(G)} d_G(v)^2 + 4f - m - 5.
\]
Equality holds if and only if $G \cong C_7$, where $f$ denotes the number of $4$-cycles in $G$.
\end{theorem}

In fact, \eqref{e1} also establishes a lower bound on the number of $4$-cycles in non-bipartite $\{C_3, C_5\}$-free graphs of size $m$. One may also notice that the bound in Theorem~\ref{thm:2.5} is attainable only when $m=7.$ Li, Peng \cite{Li-Peng} and Sun, Li \cite{sl2023}, independently, extended Theorems~\ref{thm:2.4} and \ref{thm:2.5} as follows.
\begin{theorem}[Li, Peng~\cite{Li-Peng}, Sun, Li~\cite{sl2023}]\label{thm:2.6}
Let $G$ be a non-bipartite graph in $\mathcal{G}(m, \{C_3, C_5\})$, one has $
\lambda(G) \leqslant \lambda(S_3(K_{2, \frac{m-3}{2}}))% < \sqrt{m - 3}
.
$ Equality holds if and only if $m$ is odd and $G \cong S_3(K_{2, \frac{m-3}{2}})$.
\end{theorem}

Similarly, one may see in Theorem~\ref{thm:2.6} that the corresponding extremal graph is well defined only when $m$ is odd. Sun and Li~\cite{sl2023} proposed the following problem:
%\begin{problem}[Zhai, Shu~\cite{26}, 2022]\label{pb2}
%For even $m$, what is the extremal graph among all $C_3$-free non-bipartite graphs with $m$ edges, achieving the maximum spectral radius?
%\end{problem}
\begin{problem}[Sun and Li~\cite{sl2023}, 2023]\label{pb3}
For even $m$, what is the extremal graph among all $\{C_3,C_5\}$-free non-bipartite graphs with $m$ edges, achieving the maximum spectral radius?
\end{problem}

Li and Yu~\cite{LY2023} solved Problem~\ref{pb3} for even $m\geqslant 150$.  
\begin{theorem}[Li and Yu~\cite{LY2023}, 2023]\label{thm:2.8}
Let \( m \geqslant 150 \) be even. Assume that $G$ is a non-bipartite graph in $\mathcal{G}(m, \{C_3, C_5\})$.
\begin{wst}
    \item[{\rm (a)}] If \( m \equiv 0 \pmod{3} \), then \( \lambda(G) \leqslant \lambda(S_3(K_{3, \frac{m-3}{3}})) \), with equality if and only if \( G \cong S_3(K_{3, \frac{m-3}{3}}) \).
    \item[{\rm (b)}] If \( m \equiv 1 \pmod{3} \), then \( \lambda(G) \leqslant \lambda(S_3^+(K_{2, \frac{m-4}{2}})) \), with equality if and only if \( G \cong S_3^+(K_{2, \frac{m-4}{2}}) \).
    \item[{\rm (c)}] If \( m \equiv 2 \pmod{3} \), then \( \lambda(G) \leqslant \lambda(S_3^-(K_{3, \frac{m-2}{3}})) \), with equality if and only if \( G \cong S_3^-(K_{3, \frac{m-2}{3}}) \).
\end{wst}
\end{theorem}

Zhai, Lin and Shu \cite{Zhai-Lin-Shu} showed that for $m \geqslant 8$, $K_2\vee \frac{m-1}{2}K_1$ is the unique graph among $\mathcal{G}(m, C_5)$ having the largest spectral radius (see Theorem~\ref{thm:3.5} below); the extremal graph $K_2\vee \frac{m-1}{2}K_1$ is well defined only when $m$ is odd. Min, Lou and Huang~\cite{GLH2022} extended Zhai-Lin-Shu's result mentioned above, showed that for even $m \geqslant 14$, $S_{\frac{m}{2}+2,2}^1$ is the unique graph among $\mathcal{G}(m, C_5)$ having the largest spectral radius, where $S_{\frac{m}{2}+2,2}^1$ is obtained from $K_2\vee \frac{m}{2}K_1$ by deleting an edge incident to a vertex of degree $2$ (see Theorem~\ref{thm:3.6} below). Note that $G$ is $C_5$-free if it is $\{C_5, C_7, \ldots, C_{2k+1}\}$-free. Since both $K_2\vee \frac{m-1}{2}K_1$ and $S_{\frac{m}{2}+2,2}^1$ are $\{C_5, C_7, \ldots, C_{2k+1}\}$-free, one may deduce that $K_2\vee \frac{m-1}{2}K_1$ is the unique graph among $\mathcal{G}(m, \{C_5, C_7, \ldots, C_{2k+1}\})$ achieving the maximum spectral radius, provided that $m \geqslant 8$; and for even $m,$ $S_{\frac{m}{2}+2,2}^1$ is the unique graph among $\mathcal{G}(m, \{C_5, C_7, \ldots, C_{2k+1}\})$ achieving the maximum spectral radius, provided that $m \geqslant 14$.

Note that the extremal graph in Theorem~\ref{thm:2.6} contains a $C_7.$ Lou, Lu and Huang \cite{Lou-Lu-Huang} extended Theorems~\ref{thm:2.4} and \ref{thm:2.6} to obtain the following result, which was proposed by Li and Peng \cite{Li-Peng} as a conjecture.
\begin{theorem}[Lou, Lu and Huang~\cite{Lou-Lu-Huang}, 2024]\label{thm:2.10}
For $k \geqslant 2$, if $G$ is a non-bipartite graph in $\mathcal{G}(m, \{C_3, C_5, \ldots, C_{2k+1}\})$, then
$$
\lambda(G) \leqslant \lambda(S_{2k-1}(K_{2, \frac{m-2k+1}{2}})).
$$
Equality holds if and only if $m$ is odd and $G \cong S_{2k-1}(K_{2, \frac{m-2k+1}{2}})$.
\end{theorem}

One may also see that the extremal graph in Theorem~\ref{thm:2.10} can be achieved only for odd $m$. Hence, Li, Feng and Peng \cite{Li-Feng-Peng} proposed the following question.
\begin{problem}[Li, Feng and Peng \cite{Li-Feng-Peng}, 2024]
For even $m$, what is the extremal graph attaining the maximum spectral radius over all non-bipartite $\{C_3, C_5, \ldots, C_{2k+1}\}$-free graphs with $m$ edges?
\end{problem}

By observing the extremal graphs in \cite[Theorem~7]{Li-Feng-Peng} and Theorem~\ref{thm:2.8}, Li and Yu~\cite{LY2023} proposed the following conjecture.
\begin{conjecture}[Li and Yu~\cite{LY2023}, 2023]\label{conj:2.11}
Let $k \geqslant 1$. There is an integer $N(k)$. Assume that $G$ is a non-bipartite graph in $\mathcal{G}(m, \{C_3, C_5, \ldots, C_{2k+1}\})$ with even $m \geqslant N(k)$.
\begin{wst}
    \item[{\rm (a)}] If $m - (2k - 1) \equiv 0 \pmod{3}$, then $\lambda(G) \leqslant \lambda(S_{2k-1}(K_{3, \frac{m-(2k-1)}{3}}))$, with equality if and only if $G \cong S_{2k-1}(K_{3, \frac{m-(2k-1)}{3}})$.
    \item[{\rm (b)}] If $m - (2k - 1) \equiv 1 \pmod{3}$, then $\lambda(G) \leqslant \lambda(S_{2k-1}^+(K_{2, \frac{m-(2k-1)-1}{2}}))$, with equality if and only if $G \cong S_{2k-1}^+(K_{2, \frac{m-(2k-1)-1}{2}})$.
    \item[{\rm (c)}] If $m - (2k - 1) \equiv 2 \pmod{3}$, then $\lambda(G) \leqslant \lambda(S_{2k-1}^-(K_{3, \frac{m-(2k-1)+1}{3}}))$, with equality if and only if $G \cong S_{2k-1}^-(K_{3, \frac{m-(2k-1)+1}{3}})$.
\end{wst}
\end{conjecture}
Very recently, Lou, Lu and Zhai \cite[Theorem 1.4]{LLZ2025} confirmed Conjecture~\ref{conj:2.11} with $m \geqslant 30k+45$. %Zou, Li and Feng~\cite{HFL2026} considered Problem \ref{problem1} for graphs without non-consecutive odd cycles.
%\begin{theorem}[Zou, Li, Feng~\cite{HFL2026}, 2026]
%Let \( 1 \leqslant \ell \leqslant k \) and \( G \) be a \( \{C_3, C_5, \ldots, C_{2\ell-1}, C_{2k+1}\} \)-free graph on \( n \) vertices, where \( n \geq 187k \).
%\begin{wst}
%    \item[{\rm (a)}] If \( \ell = k \) and \( \lambda(G) \geqslant \lambda(S_{2k-1}(T_{n-2k+1,2})) \), then \( G \) is bipartite, unless \( G = S_{2k-1}(T_{n-2k+1,2}) \).
%    \item[{\rm (b)}] If \( \ell \leqslant k - 1 \) and \( \lambda(G) \geqslant \lambda(C_{2\ell+1}(T_{n-2\ell,2})) \), then \( G \) is bipartite, unless \( G = C_{2\ell+1}(T_{n-2\ell,2}) \).
%\end{wst}
%\end{theorem}

\section{\normalsize Forbidding a subgraph whose chromatic number is at most 3}
In this section, we survey results on Problem~\ref{problem1} with respect to the forbidding subgraph whose chromatic number is at most 3. 

%Recall that Nosal \cite{Nosal} obtained Theorem~\ref{thm:2.1} in 1970. Over twenty years later, Nikiforov generalized Theorem~\ref{thm:2.1} to general case and characterize the corresponding extremal graphs
Nikiforov~\cite{Nikiforov3} considered Problem~\ref{problem1} for $C_4$, and established the following result.
\begin{theorem}[Nikiforov~\cite{Nikiforov3}, 2009]\label{thm4.1}
For any graph $G \in \mathcal{G}(m, C_4)$ with $m \geqslant 10$, one has $\lambda(G) \leqslant \sqrt{m}$. Equality holds if and only if $G \cong K_{1,m}$.
\end{theorem} 

Note that $K_{1,m}$ is a bipartite graph. For integers $a\geqslant 1$ and $b\geqslant 2,$ let $K_{a,b}^+$ be the graph obtained from $K_{a,b}$ by adding an edge within its part on $b$ vertices. Zhai and Shu~\cite{26} extended Theorem~\ref{thm4.1} to the following form.
\begin{theorem}[Zhai and Shu~\cite{26}, 2022]\label{thm4.2}
For any graph $G \in \mathcal{G}(m, C_4)$ with $m \geqslant 26$. If $G$ is non-bipartite, then $\lambda(G) \leqslant \lambda(K_{1,m-1}^+)$. Equality holds if and only if $G \cong K_{1,m-1}^+$.
\end{theorem}  
%For the even cycles of lengthes at least $6.$ Note that $\mathcal{G}(m, C_{2k+2})\subseteq \mathcal{G}(m, \theta_{1,2,2k-1})\cup \mathcal{G}(m, \theta_{1,2,2k})$ and that $K_k\vee(\frac{m}{k}-\frac{k-1}{2})K_1$ is $C_{2k+2}$-free. Theorems~\ref{thm:3.4} and \ref{thm:3.7} give us the following results.
%\begin{corollary}
%Let $k \geqslant 2$ and $m \geqslant 4(k^2 + 3k + 1)^2$. If $G$ is a graph in $\mathcal{G}(m, C_{2k+2})$, then  
%		$
%		\lambda(G) \leqslant \frac{k-1 + \sqrt{4m - k^2 + 1}}{2},
%		$  
%with equality if and only if $G \cong K_k\vee(\frac{m}{k}-\frac{k-1}{2})K_1$.
%\end{corollary}

On the other hand, notice that $C_4\cong K_{2,2}.$ Zhai, Lin and Shu~\cite{Zhai-Lin-Shu} extended Theorem~\ref{thm4.1} to the following form.
\begin{theorem}[Zhai, Lin and Shu~\cite{Zhai-Lin-Shu}, 2021]\label{thm4.3}
For any graph $G \in \mathcal{G}(m, K_{2,r+1})$ with $r\geqslant 2$ and $m \geqslant 16r^2,$ one has $\lambda(G) \leqslant \sqrt{m}$. Equality holds if and only if $G \cong K_{1,m}$.
\end{theorem}
Fang and You~\cite{Fang-You} extended Theorems~\ref{thm4.2} and \ref{thm4.3} by showing the following result.
\begin{theorem}[Fang and You~\cite{Fang-You}, 2023]\label{thm4.003}
For any graph $G \in \mathcal{G}(m, K_{2,r+1})\backslash \{K_{1,m}\}$ with $r\geqslant 1$ and $m \geqslant (4r+2)^2+1,$ one has $\lambda(G) \leqslant \lambda(K_{1,m-1}^+)$. Equality holds if and only if $G \cong K_{1,m-1}^+$.
\end{theorem}

Let $B_{2, r+1}=K_2\vee (r+1)K_1$ denote a book graph, which may be obtained from $K_{2,r+1}$ by adding an edge connecting two vertices of degree $r+1$. Zhai, Lin and Shu~\cite{Zhai-Lin-Shu} posed the following conjecture on book graph.
\begin{conjecture}[Zhai, Lin and Shu~\cite{Zhai-Lin-Shu}, 2021]\label{conj:3.01}
For any graph \( G\in \mathcal{G}(m, B_{2, r+1}) \) with \( m \geqslant m_0(r) \) for large enough \( m_0(r) \), one has $\lambda(G) \leqslant \sqrt{m}.$ Equality holds if and only if \( G \) is a complete bipartite graph.
\end{conjecture}
Let \( bk(G) \) stand for the \textit{booksize} of \( G \), that is, the maximum number of triangles with a common edge in \( G \).
Nikiforov~\cite{Nikiforov4} confirmed the above conjecture by showing the following result. 
\begin{theorem}[Nikiforov~\cite{Nikiforov4}, 2021]\label{thm:3.2}
If \( G \) is a graph with \( m \) edges and \( \lambda(G) \geqslant \sqrt{m} \), then
        $
        bk(G) > \frac{1}{12} \sqrt[4]{m},
        $
unless \( G \) is a complete bipartite graph.
\end{theorem}
From Theorem~\ref{thm:3.2}, one may see Conjecture~\ref{conj:3.01} holds for \( m_0(r)\geqslant (12r)^4.\) Recently, Li, Liu and Zhang \cite{Li-Liu-Zhang} improved the Nikiforov's bound on \( bk(G) \).
\begin{theorem}[Li, Liu and Zhang~\cite{Li-Liu-Zhang}, 2025]\label{thm:3.3}
If \( G \) is a graph with \( m \) edges and \( \lambda(G) \geqslant \sqrt{m} \), then
        $
        bk(G) > \frac{1}{144} \sqrt{m},
        $
unless \( G \) is a complete bipartite graph.
\end{theorem}
Theorem~\ref{thm:3.3} showed Conjecture~\ref{conj:3.01} for \( m_0(r)\geqslant (144r)^2.\) Li, Liu and Zhang~\cite{Li-Liu-Zhang} also showed the order $\sqrt{m}$ in Theorem~\ref{thm:3.3} is the best possible. 

Note that the extremal graphs in Theorem~\ref{thm:3.3} are bipartite. Motivated by Theorems~\ref{thm4.2} and \ref{thm4.003}, one may posed the following problem:
\begin{problem}\label{pb5}
What is the extremal graph attaining the maximum spectral radius over all $B_{2, r+1}$-free non-bipartite graphs with $m$ edges?
\end{problem}
Liu and Miao~\cite{LLM2025} conjectured the extremal graph in Problem~\ref{pb5} is $K_{1,m-1}^+$. But observing that $K_{r,\frac{m-1}{r}}^+$ is $B_{2, r+1}$-free, and a simple calculation shows $\lambda(K_{r,\frac{m-1}{r}}^+)>\lambda(K_{1,m-1}^+)$ provided that $m\geqslant r^2$ and $m\equiv 1 \pmod{r}$. We believe the following conjecture is true.
\begin{conjecture}%[Guo, Li, and Yu~\cite{GLY2025}, 2025]
Let $r\geqslant 1$ and $m$ be sufficiently large. If \( G \) is a \( B_{2,r+1} \)-free non-bipartite graph with size \( m \), then
$
\lambda(G) \leqslant \rho(m,r),
$ 
where $\rho(m,r)$ is the largest zero of $x^3-x^2-(m-1)x+m-2r-1=0.$
Equality holds if and only if $m\equiv 1 \pmod{r}$ and \( G \cong K_{r,\frac{m-1}{r}}^+\).
\end{conjecture}

%This result laid the groundwork for further developments, including Nikiforov’s clique-oriented extension , which provided a complete characterization of the extremal graphs. Over the years, this problem has garnered significant attention, resulting in substantial contributions from renowned researchers such as Bollob\'as, Nikiforov, and others. The extensive body of work that has emerged in this domain is well-documented in \cite{Ando,Bollo2,Edwards,Elphick,Li-Feng-Peng,Li-Peng,Lin-Ning-Wu,Lou-Lu-Huang} and their references, offering an in-depth introduction to the topic for interested readers.
For cycles of lengths $5$ and $6$, Zhai, Lin and Shu~\cite{Zhai-Lin-Shu} showed the following result.
\begin{theorem}[Zhai, Lin and Shu~\cite{Zhai-Lin-Shu}, 2021]\label{thm:3.5}
If either \( G \in \mathcal{G}(m, C_5) \) with \( m \geqslant 8 \), or \( G \in \mathcal{G}(m, C_6) \) with \( m \geqslant 22 \), then \( \lambda(G) \leqslant \frac{1+\sqrt{4m-3}}{2} \). Equality holds if and only if \( G \cong K_2\vee \frac{m-1}{2}K_1.\)
\end{theorem}
Based on Theorem~\ref{thm:3.5}, Zhai, Lin and Shu~\cite{Zhai-Lin-Shu} posed the following conjecture.
\begin{conjecture}[Zhai, Lin and Shu~\cite{Zhai-Lin-Shu}, 2021]\label{eqivZhaiconjecture}
Let $k \geqslant 2$ be a fixed integer. If $G \in \mathcal{G}(m, C_{2k+1}) \cup \mathcal{G}(m, C_{2k+2})$ and $m$ is sufficiently large, then $\lambda(G) \leqslant \frac{k-1+\sqrt{4m-k^2+1}}{2}.$ Equality holds if and only if $G \cong K_k\vee(\frac{m}{k}-\frac{k-1}{2})K_1$.
\end{conjecture}
Recall that for \( q \geqslant p \geqslant 2 \), the theta graph $\theta_{1,p,q}$ is obtained by connecting two distinct vertices with three internally disjoint paths of lengths $1,p,q$, respectively. Observe that $\theta_{1,2,2}$ is the book graph on 4 vertices. Theorem~\ref{thm:3.3} tells us that if $m\geqslant 12^4$ and $G\in \mathcal{G}(m, \theta_{1,2,2})$, then $\lambda(G) \leqslant \sqrt{m}.$ Equality holds if and only if $G$ is a complete bipartite graph. Before long, Y.T. Li advanced a stronger conjecture than Conjecture~\ref{eqivZhaiconjecture} (see also \cite{Lou-Lu-Huang}) as follows:
\begin{conjecture}[\cite{Lou-Lu-Huang}]\label{ytLiconjecture}
Let $k \geqslant 2$ be a fixed integer. If $G \in \mathcal{G}(m, \theta_{1,2,2k-1}) \cup \mathcal{G}(m, \theta_{1,2,2k})$ and $m$ is sufficiently large, then $\lambda(G) \leqslant \frac{k-1+\sqrt{4m-k^2+1}}{2}.$ Equality holds if and only if $G \cong K_k\vee(\frac{m}{k}-\frac{k-1}{2})K_1$.
\end{conjecture}
Note that $\mathcal{G}(m, C_{2k+1}) \cup \mathcal{G}(m, C_{2k+2})\subseteq \mathcal{G}(m, \theta_{1,2,2k-1}) \cup \mathcal{G}(m, \theta_{1,2,2k}).$ If Conjecture~\ref{ytLiconjecture} holds, then Conjecture~\ref{ytLiconjecture} implies Conjecture~\ref{eqivZhaiconjecture}.

As an extension of Theorem~\ref{thm:3.5}, Sun, Li and Wei~\cite{Sun-Li-Wei} solved Conjecture~\ref{ytLiconjecture} for $k=2$ by showing the following result.
\begin{theorem}[Sun, Li and Wei \cite{Sun-Li-Wei}, 2023]\label{thm:3.4}
Let \( G \) be a graph in \( \mathcal{G}(m, \theta_{1,2,r}) \). Then the following holds.
\begin{wst}
\item[{\rm (i)}] If \( r = 3 \) and \( m \geqslant 8 \), then \( \lambda(G) \leqslant \frac{1+\sqrt{4m-3}}{2} \) and equality holds if and only if \( G \cong K_2\vee \frac{m-1}{2}K_1.\)

\item[{\rm (ii)}] If \( r = 4 \) and \( m \geqslant 22 \), then \( \lambda(G) \leqslant \frac{1+\sqrt{4m-3}}{2} \) and equality holds if and only if \( G \cong K_2\vee \frac{m-1}{2}K_1.\)
\end{wst}
\end{theorem}
Lu, Lu and Li~\cite{Lu-Lu-Li} considered Conjecture~\ref{ytLiconjecture} for the case $k = 3$. %They showed $K_3\vee\frac{m-3}{3}K_1$ achieving the maximum spectral radius among $\mathcal{G}(m, \theta_{1,2,5})$ provided that $m\geqslant 38.$ 
\begin{theorem}[Lu, Lu and Li~\cite{Lu-Lu-Li}, 2024]
Let \( G \) be a graph in \( \mathcal{G}(m, \theta_{1,2,5})\). If \( m \geqslant 38 \), then
$
\lambda(G) \leqslant 1 + \sqrt{m - 2}.
$
Equality holds if and only if \(G \cong K_3 \vee \frac{m - 3}{3} K_1\).
\end{theorem}
Li, Zhai and Shu~\cite{Li-Zhai-Shu} solved Conjecture~\ref{ytLiconjecture} for $k\geqslant 3$ completely.
\begin{theorem}[Li, Zhai and Shu~\cite{Li-Zhai-Shu}, 2024]\label{thm:3.7}
Let $k \geqslant 3$ and $m \geqslant 4(k^2 + 3k + 1)^2$. If $G\in \mathcal{G}(m, \theta_{1,2,2k-1}) \cup \mathcal{G}(m, \theta_{1,2,2k})$, then  
		$
		\lambda(G) \leqslant \frac{k-1 + \sqrt{4m - k^2 + 1}}{2}.
		$  
Equality holds if and only if $G \cong K_k\vee(\frac{m}{k}-\frac{k-1}{2})K_1$.
\end{theorem}
Consequently, Conjecture \ref{eqivZhaiconjecture} was also resolved.

Furthermore, $M_{2k+2}$ and $P_{2k+2}$ are subgraphs of $\theta_{1,2,2k},$ where $M_{2k+2}$ is a graph consists of $k+1$ isolated edges. According to Theorem~\ref{thm:3.7}, Li, Zhai and Shu~\cite{Li-Zhai-Shu} obtained the following two results.
\begin{theorem}[Li, Zhai and Shu~\cite{Li-Zhai-Shu}, 2024]\label{thm:3.68}
Let $k \geqslant 3$ and $m \geqslant 4(k^2 + 3k + 1)^2$. If $G\in\mathcal{G}(m, M_{2k+2})$, then  
		$
		\lambda(G) \leqslant \frac{k-1 + \sqrt{4m - k^2 + 1}}{2}.
		$  
Equality holds if and only if $G \cong K_k\vee(\frac{m}{k}-\frac{k-1}{2})K_1$.
\end{theorem}  
\begin{theorem}[Li, Zhai and Shu~\cite{Li-Zhai-Shu}, 2024]\label{thm:3.69}
Let $k \geqslant 3$ and $m \geqslant 4(k^2 + 3k + 1)^2$. If $G\in\mathcal{G}(m, P_{2k+2})$, then  
		$
		\lambda(G) \leqslant \frac{k-1 + \sqrt{4m - k^2 + 1}}{2}.
		$  
Equality holds if and only if $G \cong K_k\vee(\frac{m}{k}-\frac{k-1}{2})K_1$.
\end{theorem}   
Wang and Wang~\cite{Wang-Wang} considered Problem~\ref{problem1} for paths of odd order.
\begin{theorem}[Wang and Wang~\cite{Wang-Wang}, 2025]\label{thm:3.70}
Let $k \geqslant 4$ and $m \geqslant 4(k+ 1)^4$. If $G\in\mathcal{G}(m, P_{2k+1})$, then  
		$
		\lambda(G) \leqslant \frac{k-2 + \sqrt{4m - k^2 + 2k}}{2}.
		$  
Equality holds if and only if $G \cong K_{k-1}\vee(\frac{m}{k-1}-\frac{k-2}{2})K_1$.
\end{theorem}

Notice that the extremal graph in Theorems~\ref{thm:3.5} and \ref{thm:3.4} is well defined only when \( m \) is odd. Let \( S_{n,k}^t \) be the graph obtained from $K_k\vee (n-t-k)K_1$ by attaching \( t \) pendant vertices to a maximum degree vertex of \(K_k\vee (n-t-k)K_1\), and \( S_n^k \) be the graph obtained from \( K_{1,n-1} \) by adding \( k \) disjoint edges within its independent set. Min, Lou and Huang~\cite{GLH2022} extended Theorem~\ref{thm:3.5} by showing the following result.%we want to obtain a sharp upper bound of \( \rho(G) \). Our result is presented as follows.
\begin{theorem}[Min, Lou and Huang~\cite{GLH2022}, 2022]\label{thm:3.6}
Let \( \tilde{p}(m) \) be the largest root of \( x^4 - mx^2 - (m - 2)x + \left( \frac{m}{2} - 1 \right) = 0 \) and \( G \) be a graph of even size \( m \) without isolated vertices. If \( G \) is a \( C_5 \)-free graph with \( m \geqslant 14 \) or \( C_6 \)-free graph with \( m \geqslant 74 \), then \( \lambda(G) \leqslant \tilde{p}(m) \), with equality if and only if \( G \cong S_{\frac{m}{2}+2,2}^1\).
\end{theorem}

Further on, Sun, Li and Wei~\cite{Sun-Li-Wei} utilized a unified approach (regardless of the parity of \( m \)) to determine all the graphs with maximum spectral radius among \( \mathfrak{G}(m, C_5) := \mathcal{G}(m, C_5) \setminus \{K_2\vee \frac{m-1}{2}K_1\} \) or \( \mathfrak{G}(m, C_6) := \mathcal{G}(m, C_6) \setminus \{K_2\vee \frac{m-1}{2}K_1\} \) for \( m \geqslant 22 \).

Let \( \rho_1(m) \) be the largest zero of \( \psi_1(x) \), where
$$
\psi_1(x) = 
\begin{cases} 
x^4 - mx^2 - (m - 2)x + \frac{m}{2} - 1, & \text{if } m \text{ is even}; \\
x^4 - mx^2 - (m - 3)x + m - 3, & \text{if } m \text{ is odd},
\end{cases}
$$
and let \( \rho_2(m) \) be the largest zero of \( \psi_2(x) \), where
$$
\psi_2(x) = 
\begin{cases} 
x^3 - 2x^2 - (m - 3)x + m - 6, & \text{if } m \text{ is even and } 22 \leqslant m \leqslant 72; \\
x^4 - mx^2 - (m - 2)x + \frac{m}{2} - 1, & \text{if } m \text{ is even and } m \geqslant 74; \\
x^5 - x^4 - (m - 1)x^3 - 2x^2 + \frac{3m - 17}{2}x - \frac{m - 7}{2}, & \text{if } m \text{ is odd and } 23 \leqslant m \leqslant 71; \\
x^4 - mx^2 - (m - 3)x + m - 3, & \text{if } m \text{ is odd and } m \geqslant 73.
\end{cases}
$$

The subsequent two results identify all the graphs having the largest spectral radii among \( \mathfrak{G}(m, C_5) \) and \( \mathfrak{G}(m, C_6) \), respectively.
\begin{theorem}[Sun, Li and Wei \cite{Sun-Li-Wei}, 2023]\label{thm:3.04}
Let \( G \) be in \( \mathfrak{G}(m, C_5) \) with \( m \geqslant 22 \). Then \( \lambda(G) \leqslant \rho_1(m) \). Equality holds if and only if \( G \cong S_{\frac{m+4}{2}, 2}^1 \) if \( m \) is even and \( G \cong S_{\frac{m+5}{2}, 2}^2 \) if \( m \) is odd.
\end{theorem}
\begin{theorem}[Sun, Li and Wei \cite{Sun-Li-Wei}, 2023]\label{thm:3.06}
Let \( G \) be in \( \mathfrak{G}(m, C_6) \) with \( m \geqslant 22 \). Then \( \lambda(G) \leqslant \rho_2(m) \). Equality holds if and only if
$$
G \cong 
\begin{cases} 
K_1 \vee S_{\frac{m}{2}}^1, & \text{if } m \text{ is even and } 22 \leqslant m \leqslant 72; \\
S_{\frac{m+4}{2}, 2}^1, & \text{if } m \text{ is even and } m \geqslant 74; \\
K_1 \vee (S_{\frac{m-1}{2}}^1 \cup K_1), & \text{if } m \text{ is odd and } 23 \leqslant m \leqslant 71; \\
S_{\frac{m+5}{2}, 2}^2, & \text{if } m \text{ is odd and } m \geqslant 73.
\end{cases}
$$
\end{theorem}

Fang, You~\cite{Fang-You} and Liu, Wang~\cite{Liu-Wang} extended Theorems~\ref{thm:3.4}, \ref{thm:3.6}, \ref{thm:3.04} and \ref{thm:3.06} by showing the following results.
\begin{theorem}[Fang and You \cite{Fang-You}, 2023]\label{thm:3.004}
Let \( G \) be in \( \mathcal{G}(m, \theta_{1,2,3})\setminus \{K_2\vee \frac{m-1}{2}K_1\} \) with \( m \geqslant 22 \). Then \( \lambda(G) \leqslant \rho_1(m) \). Equality holds if and only if \( G \cong S_{\frac{m+4}{2}, 2}^1 \) if \( m \) is even and \( G \cong S_{\frac{m+5}{2}, 2}^2 \) if \( m \) is odd.
\end{theorem}
\begin{theorem}[Liu and Wang \cite{Liu-Wang}, 2024]\label{thm:3.006}
Let \( G \) be in \( \mathcal{G}(m, \theta_{1,2,4}) \) with even \( m \geqslant 74 \). Then \( \lambda(G) \leqslant \rho_1(m) \). Equality holds if and only if \( G \cong S_{\frac{m+4}{2}, 2}^1 \).
\end{theorem}
It is worth noting that the extremal graph $K_k\vee(\frac{m}{k}-\frac{k-1}{2})K_1$ in Theorem~\ref{thm:3.7} is well-defined only when $m + \frac{k(k+1)}{2} \equiv 0 \pmod{k}$. It is natural to raise the following problem. 
\begin{problem}[Liu et al.~\cite{LLLY2025}, 2026]\label{pb7}
For $k\geqslant 3,$ when $m + \frac{k(k+1)}{2} \equiv l \pmod{k}$ with $1 \leqslant l \leqslant k-1$, which graphs in $\mathcal{G}(m, \theta_{1,2,2k-1}) \cup \mathcal{G}(m, \theta_{1,2,2k}),$ achieve the largest spectral radius?   
\end{problem}
Let $S^{1,l}_{n,k} := K_l\vee((K_{k-l}\vee (n-k-1)K_1)\cup K_1).$
By employing $k$-core analysis combined with spectral techniques, Liu et al.~\cite{LLLY2025} overcomed Problem~\ref{pb7} as follows.
\begin{theorem}[Liu et al.~\cite{LLLY2025}, 2026]\label{mainresult}
Assume $1 \leqslant l < k$ and $m \geqslant 144k^4$ with $k \geqslant 3$. If $G\in\mathcal{G}(m, \theta_{1,2,2k-1}) \cup \mathcal{G}(m, \theta_{1,2,2k}),$ then $\lambda(G) \leqslant \lambda(S^{1,l}_{\frac{m-l}{k} + \frac{k+3}{2}, k})$ for $m+\frac{k(k+1)}{2}\equiv l \pmod{k}$. Equality holds if and only if $G \cong S^{1,l}_{\frac{m-l}{k} + \frac{k+3}{2}, k}$.
\end{theorem}
Consequently, Liu et al.~\cite{LLLY2025} derived the following result.
\begin{theorem}[Liu et al.~\cite{LLLY2025}, 2026]\label{mainresult2}
Assume $1 \leqslant l < k$ and $m \geqslant 144k^4$ with $k \geqslant 3$. If $G\in\mathcal{G}(m, C_{2k+1}) \cup \mathcal{G}(m, C_{2k+2})$, then  
		$
		\lambda(G) \leqslant \lambda(S^{1,l}_{\frac{m-l}{k} + \frac{k+3}{2}, k})
		$  
for $m+\frac{k(k+1)}{2}\equiv l\pmod{k}$. Equality holds if and only if $G \cong S^{1,l}_{\frac{m-l}{k} + \frac{k+3}{2}, k}$.
\end{theorem}
Similarly to Theorems~\ref{thm:3.68} and \ref{thm:3.69}, we may derive the following two results from Theorems~\ref{thm:3.004}, \ref{thm:3.006} and \ref{mainresult}.
\begin{corollary}
Assume $1 \leqslant l < k$ and $m \geqslant 144k^4$ with $k \geqslant 2$. If $G\in\mathcal{G}(m, M_{2k+2}),$ then $\lambda(G) \leqslant \lambda(S^{1,l}_{\frac{m-l}{k} + \frac{k+3}{2}, k})$ for $m+\frac{k(k+1)}{2}\equiv l \pmod{k}$. Equality holds if and only if $G \cong S^{1,l}_{\frac{m-l}{k} + \frac{k+3}{2}, k}$.
\end{corollary}
\begin{corollary}
Assume $1 \leqslant l < k$ and $m \geqslant 144k^4$ with $k \geqslant 2$. If $G\in\mathcal{G}(m, P_{2k+2}),$ then $\lambda(G) \leqslant \lambda(S^{1,l}_{\frac{m-l}{k} + \frac{k+3}{2}, k})$ for $m+\frac{k(k+1)}{2}\equiv l \pmod{k}$. Equality holds if and only if $G \cong S^{1,l}_{\frac{m-l}{k} + \frac{k+3}{2}, k}$.
\end{corollary}
%We may also observe that the extremal graph $K_{k-1}\vee(\frac{m}{k-1}-\frac{k-2}{2})K_1$ in Theorem~\ref{thm:3.70} is well-defined only when $m + \frac{k(k-1)}{2} \equiv 0 \pmod{k-1}$. Motivated by the structures of the extremal graphs in both Theorem~\ref{thm:3.70} and  Theorem~\ref{mainresult}, we believe the following conjecture is true.
%\begin{conjecture}
%Assume $1 \leqslant l < k-1$ with $k \geqslant 4$. There is an $m_0$ such that for all $m\geqslant m_0,$ if $G\in\mathcal{G}(m, P_{2k+1})$, then $\lambda(G) \leqslant \lambda(S^{1,l}_{\frac{m-l}{k-1} + \frac{k+2}{2}, k-1})$ for $m+\frac{k(k-1)}{2}\equiv l \pmod{k-1}$. Equality holds if and only if $G \cong S^{1,l}_{\frac{m-l}{k-1} + \frac{k+2}{2}, k-1}$. 
%\end{conjecture}

Inspired by Theorem~\ref{thm:3.7}, Li, Zhao and Zou~\cite{Li-Zhao-Zou} posed the following problem.
\begin{problem}[Li, Zhao and Zou~\cite{Li-Zhao-Zou}, 2025]\label{p6}
What is the maximum spectral radius of graphs among \( \mathcal{G}(m, \theta_{1,p,q}) \) for \( q \geqslant p \geqslant 3 \)?
\end{problem}
Gao and Li~\cite{GaoL2025} solved Problem~\ref{p6} for \( q=p=3 \).
\begin{theorem}[Gao and Li~\cite{GaoL2025}, 2025]\label{Gaoli}
Let $G\in\mathcal{G}(m, \theta_{1,3,3})$ with $m\geqslant 43.$ Then $\lambda(G)\leqslant \frac{1+\sqrt{4m-3}}{2}$, and equality holds if and only if $G\cong K_{2}\vee\frac{m-1}{2}K_{1}.$
\end{theorem}
For $k\geqslant 2,$ the friendship graph $F_{k,3}$ is defined as $F_{k,3}=K_1\vee kK_2.$ Li, Lu and Peng~\cite{LLP2023} proposed the following conjecture for graphs without $F_{k,3}$.
\begin{conjecture}[Li, Lu and Peng~\cite{LLP2023}, 2023]\label{conj:2023}
Let $k \geqslant 2$ be fixed and $m$ be large enough. If $G\in\mathcal{G}(m, F_{k,3})$, then 
$
\lambda(G) \leqslant \frac{k - 1 + \sqrt{4m - k^2 + 1}}{2}.
$ 
Equality holds if and only if $G\cong K_k\vee(\frac{m}{k}-\frac{k-1}{2})K_1$.
\end{conjecture}
Li, Lu and Peng~\cite{LLP2023} solved Conjecture~\ref{conj:2023} for $k=2.$
\begin{theorem}[Li, Lu and Peng~\cite{LLP2023}, 2023]\label{thm:3.0012}
If $G\in\mathcal{G}(m, F_{2,3})$ with $m\geqslant 8,$ then 
$
\lambda(G) \leqslant \frac{1 + \sqrt{4m-3}}{2}.
$ 
Equality holds if and only if $G\cong K_2\vee\frac{m-1}{2}K_1$.
\end{theorem}
Yu, Li and Peng~\cite{YLP2025} solved Conjecture~\ref{conj:2023} for $k=3.$
\begin{theorem}[Yu, Li and Peng~\cite{YLP2025}, 2025]
If $G\in\mathcal{G}(m,F_{3,3})$ with $m \geqslant 33$, then
$
\lambda(G) \leqslant 1 + \sqrt{m - 2}.
$
Equality holds if and only if $G\cong K_3 \vee \frac{m - 3}{3}K_1$.
\end{theorem}
Define $F_t=K_1 \vee P_{t-1}$ to be the \textit{fan graph} on $t$ vertices. It is worth noting that $F_3$ is a triangle and $F_4$ is a book graph on 4 vertices. The spectral extremal problem on $F_t\,(3\leqslant t\leqslant4)$-free graphs with given size $m$ were studied by Nosal \cite{Nosal}, Nikiforov \cite{VNK2009,Nikiforov4} and Li, Liu and Zhang~\cite{Li-Liu-Zhang}, see Theorems~\ref{thm:2.1}, \ref{thm:2.2}, \ref{thm:3.2} and \ref{thm:3.3}. 
Zhang and Wang~\cite{ZW2024} posed the following problem.
\begin{problem}[Zhang and Wang~\cite{ZW2024}]\label{pb6}
What is the maximum spectral radius and what are the corresponding extremal graphs among all $F_t$-free graphs with size $m$ for $t\geqslant 5$?
\end{problem}
Zhang, Wang~\cite{ZW2024} and Yu, Li, Peng~\cite{YLP2025} solved Problem~\ref{pb6} for $t=5$ independently. 
\begin{theorem}[Yu, Li, Peng~\cite{YLP2025} and Zhang, Wang~\cite{ZW2024}]\label{thm:3.013}
If $G\in\mathcal{G}(m,F_5)$ with $m \geqslant 8$, then
$
\lambda(G) \leqslant \frac{1 + \sqrt{4m - 3}}{2}.
$
Equality holds if and only if $G\cong K_2 \vee \frac{m - 1}{2}K_1$.
\end{theorem} 
Gao and Li~\cite{GL2025} solved Problem~\ref{pb6} for $t=6$.
\begin{theorem}[Gao and Li~\cite{GL2025}, 2026]\label{thm:3.13}
If \( G \in  \mathcal{G}(m, F_6) \) with \( m \geqslant 88 \), then \( \lambda(G) \leqslant \frac{1+\sqrt{4m-3}}{2} \). Equality holds if and only if \( G \cong K_2 \vee \frac{m-1}{2}K_1 \).
\end{theorem}
Zhang and Wang~\cite{ZW2025} solved Problem~\ref{pb6} for $t=7$.
\begin{theorem}[Zhang and Wang~\cite{ZW2025}, 2025]
If \( G \in  \mathcal{G}(m, F_7) \) with \( m \geqslant 33 \), then \( \lambda(G) \leqslant 1+\sqrt{m-2}\). Equality holds if and only if \( G \cong K_3 \vee \frac{m-3}{3}K_1 \).
\end{theorem}
Yu, Li and Peng~\cite{YLP2025} gave a conjecture on Problem~\ref{pb6} as follows.
\begin{conjecture}[Yu, Li and Peng~\cite{YLP2025}, 2025]\label{conj-6}
Let $k \geqslant 2$ be ﬁxed and $m$ be suﬃciently large. If $G\in  \mathcal{G}(m, F_{2k+1})$ or $ G\in\mathcal{G}(m, F_{2k+2}),$ then
$
\lambda(G)\leqslant\frac{k-1+\sqrt{4m-k^2+1}}{2}.
$
Equality holds if and only if $G \cong K_k\vee(\frac{m}{k}-\frac{k-1}{2})K_1$.
\end{conjecture}
Recently, Li, Zhao and Zou~\cite{Li-Zhao-Zou} solved Conjecture~\ref{conj-6} for $k\geqslant 3$ completely.
\begin{theorem}[Li, Zhao and Zou~\cite{Li-Zhao-Zou}, 2025]\label{thm3.27}
Let $k\geqslant3$ and $m\geqslant \frac{9}{4}k^6+6k^5+46k^4+56k^3+196k^2$. If $G\in\mathcal{G}(m, F_{2k+2}),$  then
$
\lambda(G)\leqslant\frac{k-1+\sqrt{4m-k^2+1}}{2}.
$
Equality holds if and only if $G \cong K_k\vee(\frac{m}{k}-\frac{k-1}{2})K_1$.
\end{theorem}
Note that $\mathcal{G}(m,F_{2k+1})\subseteq \mathcal{G}(m,F_{2k+2}),$ and $K_{k}\vee (\frac{m}{k}-\frac{k-1}{2})K_1\in\mathcal{G}(m,F_{2k+1})$. Hence, Conjecture~\ref{conj-6} follows directly from Theorems~\ref{thm:3.13} and \ref{thm3.27}.

One may see that for $q\geqslant p\geqslant 3, s\geqslant r\geqslant 3, p+q= 2k+ 1$ and $r+ s= 2k+ 2$, $\mathcal{G}( m, \theta _{1, p, q}) \cup \mathcal{G} ( m, \theta _{1, s, t}) \subseteq \mathcal{G} ( m, F_{2k+ 2})$, and there is no $\theta_{1,p,q}$ or $\theta_{1,s,t}$ in $K_k\vee(\frac{m}{k}-\frac{k-1}{2})K_1.$ By Theorem~\ref{thm3.27}, Li, Zhao and Zou~\cite{Li-Zhao-Zou} obtained the following result.
\begin{theorem}[Li, Zhao and Zou~\cite{Li-Zhao-Zou}, 2025]\label{thm3.014}
Let $k\geqslant 3$ and $m\geqslant \frac{9}{4}k^6+6k^5+46k^4+56k^3+196k^2$. If $G\in\mathcal{G}(m, \theta_{1,p,q})$ or $G\in \mathcal{G}(m, \theta_{1,r,s})$ with $q\geqslant p\geqslant 3,s\geqslant r\geqslant 3, p+q=2k+1$ and $r+s=2k+2,$ then
$
\lambda(G)\leqslant\frac{k-1+\sqrt{4m-k^2+1}}{2}.
$
Equality holds if and only if $G \cong K_k\vee(\frac{m}{k}-\frac{k-1}{2})K_1$.
\end{theorem}
Theorem~\ref{thm3.014} together with Theorem~\ref{Gaoli} solve Problem~\ref{p6} completely.

Inspired by Theorem~\ref{thm3.014}, Li, Zhao and Zou~\cite{Li-Zhao-Zou} proposed the following problem.
\begin{problem}[Li, Zhao and Zou~\cite{Li-Zhao-Zou}, 2025]\label{qu2} 
How can we characterize the graphs among \( \mathcal{G}(m, \theta_{r,p,q}) \) having the largest spectral radius for \( q \geqslant p \geqslant r \geqslant 2 \)?
\end{problem}
Note that $\theta_{2,2,2}\cong K_{2,3}.$ Zhai, Lin and Shu~\cite{Zhai-Lin-Shu} solved Problem~\ref{qu2} for \( r = p = q = 2 \) (see Theorem~\ref{thm4.3}). Recently, Gao and Li~\cite{GLAx2025} solved Problem~\ref{qu2} for \( r = 2, p = 2 \) and \( q = 3 \).
\begin{theorem}[Gao and Li \cite{GLAx2025}, 2025]\label{thm:2.40}
Let $G\in\mathcal{G}(m, \theta_{2,2,3})$ with $m\geqslant 57,$ then  
		$
		\lambda(G) \leqslant \frac{1 + \sqrt{4m - 3}}{2}.
		$  
Equality holds if and only if $G \cong K_2\vee\frac{m-1}{2}K_1$.
\end{theorem} 

One may also see that $\mathcal{G}(m,F_{k,3})\subseteq\mathcal{G}(m,F_{2k+2}),$ and there is no $F_{k,3}$ in $K_k\vee(\frac{m}{k}-\frac{k-1}{2})K_1.$ By Theorem~\ref{thm3.27}, Li, Zhao and Zou~\cite{Li-Zhao-Zou} obtained the following result.
\begin{theorem}[Li, Zhao and Zou~\cite{Li-Zhao-Zou}, 2025]\label{thm:3.50}
Let $k\geqslant3$ and $m\geqslant \frac{9}{4}k^6+6k^5+46k^4+56k^3+196k^2$. 
If $G\in\mathcal{G}(m,F_{k,3}),$ then
$
\lambda(G)\leqslant\frac{k-1+\sqrt{4m-k^2+1}}{2}.
$
Equality holds if and only if $G\cong K_k\vee(\frac{m}{k}-\frac{k-1}{2})K_1$.
\end{theorem}
Together with Theorems~\ref{thm:3.0012} and \ref{thm:3.50}, Conjecture~\ref{conj:2023} is confirmed.

Observe that the bound in Theorems~\ref{Gaoli}, \ref{thm:3.0012}, \ref{thm:3.013} and \ref{thm:3.13} is attainable only when $m$ is odd; and the bound in Theorems~\ref{thm3.27}, \ref{thm3.014} and \ref{thm:3.50} is attainable only when $m+\frac{k(k+1)}{2}\equiv 0 \pmod{k}.$ Motivated by Theorem~\ref{mainresult}, we may posed the following three conjectures.
\begin{conjecture}\label{conj:3.20}%[Li, Zhao, and Zou~\cite{Li-Zhao-Zou}, 2025]\label{thm3.014}
Let $k\geqslant 2.$ There is an $m_0$ such that for all $m\geqslant m_0,$ and for all $G\in\mathcal{G}(m, \theta_{1,p,q})\cup\mathcal{G}(m, \theta_{1,r,s})$ with $q\geqslant p\geqslant 3,s\geqslant r\geqslant 3, p+q=2k+1$ and $r+s=2k+2.$ If $m+\frac{k(k+1)}{2}\equiv l \pmod{k}$ with $1\leqslant l\leqslant k-1,$ then $\lambda(G) \leqslant \lambda(S^{1,l}_{\frac{m-l}{k} + \frac{k+3}{2}, k}).$ Equality holds if and only if $G \cong S^{1,l}_{\frac{m-l}{k} + \frac{k+3}{2}, k}$.
\end{conjecture}

\begin{conjecture}\label{conj:3.21}%[Li, Zhao, and Zou~\cite{Li-Zhao-Zou}, 2025]\label{thm3.014}
Let $k\geqslant 2.$ There is an $m_0$ such that for all $m\geqslant m_0,$ and for all $G\in\mathcal{G}(m, F_{k,3})$. If $m+\frac{k(k+1)}{2}\equiv l \pmod{k}$ with $1\leqslant l\leqslant k-1,$ then $\lambda(G) \leqslant \lambda(S^{1,l}_{\frac{m-l}{k} + \frac{k+3}{2}, k}).$ Equality holds if and only if $G \cong S^{1,l}_{\frac{m-l}{k} + \frac{k+3}{2}, k}$.
\end{conjecture}

\begin{conjecture}\label{conj:3.22}%[Li, Zhao, and Zou~\cite{Li-Zhao-Zou}, 2025]\label{thm3.014}
Let $k\geqslant 2.$ There is an $m_0$ such that for all $m\geqslant m_0,$ and for all $G\in\mathcal{G}(m, F_{2k+2})$. If $m+\frac{k(k+1)}{2}\equiv l \pmod{k}$ with $1\leqslant l\leqslant k-1,$ then $\lambda(G) \leqslant \lambda(S^{1,l}_{\frac{m-l}{k} + \frac{k+3}{2}, k}).$ Equality holds if and only if $G \cong S^{1,l}_{\frac{m-l}{k} + \frac{k+3}{2}, k}$.
\end{conjecture}
Chen and Yuan~\cite{CY2025} obtained the following result.%Conjecture~\ref{conj:3.22} for $F_5.$
\begin{theorem}[Chen and Yuan~\cite{CY2025}, 2025]\label{thm:3.60}
If \( m \geqslant 92 \) is even, and $G\in\mathcal{G}(m,F_5)$, then
$
\lambda(G) \leqslant \lambda(S^{1,1}_{\frac{m+4}{2},2}).
$
Equality holds if and only if $G \cong S^{1,1}_{\frac{m+4}{2},2}$.
\end{theorem}

Note that $F_{2,3}$ is a subgraph of $F_5,$ Theorem~\ref{thm:3.60} deduces Conjecture~\ref{conj:3.21} for $k=2$. Independently, Chen et al. \cite{CLZZ2025} confirmed Conjecture~\ref{conj:3.21} for $k=2$ as follows.
\begin{theorem}[Chen et al. \cite{CLZZ2025}, 2025]
If \( m \geqslant 16 \) is even, and $G\in\mathcal{G}(m,F_{2,3})$, then
$
\lambda(G) \leqslant \lambda(S^{1,1}_{\frac{m+4}{2},2}).
$
Equality holds if and only if $G \cong S^{1,1}_{\frac{m+4}{2},2}$.
\end{theorem}

The wheel graph on $k$ vertices is defined as $W_k=K_1\vee C_{k-1}$. We call a wheel graph odd (resp. even) if its order is odd (resp. even). Yu, Li and Peng~\cite{YLP2025} posed the following conjecture for odd wheels. 
\begin{conjecture}\label{conj:3.2}%[Yu, Li and Peng~\cite{YLP2025}, 2025]
Let \( k \geqslant 2 \) be fixed and \( m \) be sufficiently large, and \( G \) is a \( W_{2k+1} \)-free graph with \( m \) edges. 
\begin{itemize}
  \item[{\rm (i)}]If $m+\frac{k(k+1)}{2}\equiv 0 \pmod{k},$ then $\lambda(G) \leq \frac{k - 1 + \sqrt{4m - k^2 + 1}}{2},$ equality holds if and only if \( G = K_k \vee \left( \frac{m}{k} - \frac{k - 1}{2} \right) K_1 \). \hfill{(Lu, Li, Peng~{\rm \cite{YLP2025}, 2025})}
  \item [{\rm (ii)}]If $m+\frac{k(k+1)}{2}\equiv l \pmod{k}$ with $1\leqslant l\leqslant k-1,$ then $\lambda(G) \leqslant \lambda(S^{1,l}_{\frac{m-l}{k} + \frac{k+3}{2}, k}),$ equality holds if and only if $G \cong S^{1,l}_{\frac{m-l}{k} + \frac{k+3}{2}, k}$. %\hfill{(Yu \& Zhang, (2025+))}
\end{itemize}
\end{conjecture}
Conjecture~\ref{conj:3.2} is still open. Li, Liu and Zhai~\cite{LiLZ2025} obtained the following result.
\begin{theorem}[Li, Liu and Zhai~\cite{LiLZ2025}, 2025]\label{thm:3.61}
Let $G$ be a graph of size $m\geqslant 25$ with no isolated vertices. If $\lambda(G)\geqslant \frac{1 + \sqrt{4m - 5}}{2},$ then $G$ contains a wheel unless $G$ is isomorphic to either $K_2\vee\frac{m-1}{2}K_1$ or $S_{\frac{m+4}{2}, 2}^1.$ 
\end{theorem}
Motivated by Theorem~\ref{thm:3.61}, we may posed a weaker conjecture than Conjecture~\ref{conj:3.2}.
\begin{conjecture}
Let $k \geqslant 2$ be fixed and $m$ be sufficiently large. If $G\in \mathcal{G}(m, \{W_t:t\geqslant 2k+1\})$, then
$
    \lambda(G) \leqslant \frac{k - 1 + \sqrt{4m - k^2 + 1}}{2}.
$
Equality holds if and only if $G\cong K_k \vee ( \frac{m}{k} - \frac{k - 1}{2}) K_1$.
\end{conjecture}

A graph $F$ is called \textit{almost-bipartite} if it is color-critical with $\chi(F) = 3$ or bipartite.
\begin{definition}
For an almost-bipartite graph \( F \), let \(\mathcal{A}_F\) be the family of subgraphs induced by $V(F) \setminus I,$ where \( I \) ranges over the maximal independent sets of \( F \). Let \(\mathcal{M}_F\) be the family of \(\mathcal{A}_F\)-free graphs.
\end{definition}
Very recently, Li, Liu and Zhang~\cite{LiLiuZhang2} obtained the following nice result.
\begin{theorem}[Li, Liu and Zhang~\cite{LiLiuZhang2}, 2025]\label{thm:3.62}
For any almost-bipartite graph \( F \) that is not a star, if \( m \) is sufficiently large and \( G\) attains the maximum spectral radius among all graphs in \(\mathcal{G}(m, F)\), then one of the following holds:
\begin{wst}
\item[{\rm (i)}] \( G \) is a complete bipartite graph.
\item[{\rm (ii)}] {\( G \) has a vertex partition \( A \cup C \) such that \( |A| < |F| \) and \( G[A] \) is a non-empty graph in \( \mathcal{M}_F \), \( C \) is an independent set, and all but at most one vertex of \( C \) are adjacent to all vertices of \( A \).}
\end{wst}
\end{theorem}
As products of Theorem~\ref{thm:3.62}, Li, Liu and Zhang~\cite{LiLiuZhang2} solved Problem~\ref{problem1} for a large number of almost-bipartite graphs (see \cite[Table~1]{LiLiuZhang2}). In particular, Li, Liu and Zhang~\cite{LiLiuZhang2} solved Problem~\ref{pb7} with a quite different method from that of Liu et al.~\cite{LLLY2025}, and solved Problem~\ref{qu2}, Conjecture~\ref{conj:3.20} completely.

A planar graph $G$ is a graph that can be drawn on the plane such that the edges of $G$ intersect only at their endpoints. A graph $G$ is called $k$-planar if it has a drawing in the plane $\mathbb{R}^2$ such that each edge of $G$ is crossed with at most $k$ other edges. Fan, Kang and Wu~\cite{FKW} determined the largest spectral radius of a planar graph and characterized the corresponding extremal graphs. 
\begin{theorem}[Fan, Kang and Wu~\cite{FKW}, 2026]\label{thm:5.04}
Let $G$ be a planar graph with $m$ edges. Suppose that $m$ is sufficiently large, then for odd $m$, $\lambda(G)\leqslant\lambda(K_2\vee \frac{m-1}{2}K_1)$, with equality if and only if \( G \cong K_2\vee \frac{m-1}{2}K_1 \); for even $m,$ \( \lambda(G) \leqslant\lambda(S_{\frac{m+4}{2}, 2}^1) \), with equality if and only if \( G \cong S_{\frac{m+4}{2}, 2}^1 \).
\end{theorem}
It is well-known that if $G$ is planar, then $G$ is $K_{3,3}$-free; and there is a $t_0$ such that any $k$-planar graph $G$ is $K_{3,t}$-free for every $t\geqslant t_0$ (see~\cite{K1970}). As products of Theorem~\ref{thm:3.62}, Li, Liu and Zhang~\cite{LiLiuZhang2} also showed Theorem~\ref{thm:5.04} and gave the following result. 
\begin{theorem}[Li, Liu and Zhang~\cite{LiLiuZhang2}, 2026]\label{thm:5.05}
Let $G$ be a $k$-planar graph with $m$ edges. Suppose that $m$ is sufficiently large, then for odd $m$, $\lambda(G)\leqslant\lambda(K_2\vee \frac{m-1}{2}K_1)$, with equality if and only if \( G \cong K_2\vee \frac{m-1}{2}K_1 \); for even $m,$ \( \lambda(G) \leqslant\lambda(S_{\frac{m+4}{2}, 2}^1) \), with equality if and only if \( G \cong S_{\frac{m+4}{2}, 2}^1 \).
\end{theorem}

\section{\normalsize Forbidding a color-critical subgraph whose chromatic number is at least 4}
	
Nikiforov~\cite{Nikiforov1,VNK2009} solved Problem~\ref{problem1} for complete graphs.
\begin{theorem}[Nikiforov~\cite{Nikiforov1,VNK2009}]\label{thm:5.1}
For $r\geqslant 3$, let $G\in\mathcal{G}(m, K_{r+1})$. Then one has $\lambda(G) \leqslant \sqrt{2m(r-1)/r}.$ Equality holds if and only if $G$ is a complete regular $r$-partite graph.
\end{theorem} 
Bollob\'as and Nikiforov~\cite{Bollo2} further posed the following conjecture.
\begin{conjecture}[Bollob\'as and Nikiforov~\cite{Bollo2}]\label{conj:5.1}
Let $G$ be a $K_{r+1}$-free graph of size $m\geqslant r + 1.$ Then 
$$
\lambda_1^2(G)+\lambda_2^2(G)\leqslant 2m(1-\frac{1}{r}).
$$
\end{conjecture} 
Lin, Ning and Wu~\cite{9} solved Conjecture~\ref{conj:5.1} for $r=2$ (see Theorem~\ref{thm:5.2}), and Conjecture~\ref{conj:5.1} is still open for $r\geqslant 3.$

Note that the graph in $\mathcal{G}(m, K_{r+1})$ with maximum spectral radius is $r$-partite. Lou, Lu and Zhai \cite{LLZ2025} posed the following problem. 
\begin{problem}[Lou, Lu and Zhai~\cite{LLZ2025}, 2025]\label{pb:5.3} 
For every $r \geqslant 3$ and every positive integer $m = ar + b$, where $a$ is a non-negative integer and $1 \leqslant b \leqslant r$, what is the extremal graph with maximal spectral radius over all non-$r$-partite $K_{r+1}$-free graphs on $m$ edges?
\end{problem}
For $k\geqslant 1,\,r\geqslant 3,$ the generalized book graph is defined as \( B_{r,k+1}= K_r \vee (k+1)K_1.\) Clearly, \( B_{r,k+1} \) contains \( k+1 \) copies of \( K_{r+1} \). Li, Liu and Feng~\cite{Li-Liu-Feng} posed a more general conjecture than Theorem~\ref{thm:5.1}.
\begin{conjecture}[Li, Liu and Feng~\cite{Li-Liu-Feng}, 2022]\label{conj:5.2} 
For $k\geqslant 1,\,r\geqslant 3$ and sufficiently large $m.$ If \( G \in \mathcal{G}(m, B_{r,k+1}) \), then $\lambda(G) \leqslant \sqrt{2m(r-1)/r}.$ Equality holds if and only if \( G \) is a complete regular \( r \)-partite graph.
\end{conjecture}
Li, Liu and Zhang~\cite{Li-Liu-Zhang} solved Conjecture~\ref{conj:5.2} by showing the following result.
\begin{theorem}[Li, Liu and Zhang~\cite{Li-Liu-Zhang}, 2025]\label{thm:5.02}
Every \( m \)-edge graph \( G \) with \( \lambda^2(G) > (1 - \frac{1}{r})\cdot 2m\) contains a copy of \( B_{r,k+1} \) with \( k = \Omega_{r}(\sqrt{m}) \). Furthermore, there are such graphs with largest generalized book size \( O_{r}(\sqrt{m}) \).
\end{theorem} 
For even wheels, Yu, Li and Peng~\cite{YLP2025} posed the following conjecture.
\begin{conjecture}[Yu, Li and Peng~\cite{YLP2025}, 2025]\label{conj:3.3}
Let $k \geqslant 2$ be fixed and $m$ be large enough. If $G\in \mathcal{G}(m, W_{2k+2})$, then
$
    \lambda(G) \leqslant \sqrt{4m/3}.
$
Equality holds if and only if $G$ is a regular complete $3$-partite graph.
\end{conjecture}

Recently, Li, Liu and Zhang~\cite{LiLiuZhang} gave an edge-spectral version of Erd\"os-Stone-Simonovits theorem, resolved Problem~\ref{problem1} asymptotically for any graph $F$.
\begin{theorem}[Li, Liu and Zhang~\cite{LiLiuZhang}, 2025]\label{thm:5.03}
If \(\chi(F) = r + 1 \geqslant 3\) and \( G \in \mathcal{G}(m, F) \), then
$$
\lambda^2(G) \leqslant (1 - \frac{1}{r} + o(1))\cdot 2m.
$$
\end{theorem} 
Note that all of $K_{r+1},\,B_{r,k+1}$ and $W_{2k+2}$ are color-critical. Motivated by Theorem~\ref{thm:5.1}, Conjectures~\ref{conj:5.2} and \ref{conj:3.3}, Yu and Li~\cite{YL2025} posed the following problem.
\begin{problem}\label{pb:5.2}
Given a color-critical graph \( F \) with \( \chi(F)=r + 1 \), determine the graphs that attain the maximum spectral radius among all \( F \)-free graphs with \( m \) edges.
\end{problem} 
Recently, Li, Liu and Zhang~\cite{LiLiuZhang2} solved Problem~\ref{pb:5.2} by showing the following nice result.
\begin{theorem}[Li, Liu and Zhang~\cite{LiLiuZhang2}, 2025]\label{thm:5.03}
Let \( F \) be a color-critical graph with \( \chi(F) = r + 1 \geqslant 4 \). For sufficiently large \( m \), if \( G \in \mathcal{G}(m, F )\), then
$$
\lambda_1^2(G) \leqslant ( 1 - \frac{1}{r} ) 2m.
$$
Equality holds if and only if \( G \) is a regular complete \( r \)-partite graph.
\end{theorem}
Consequently, Conjectures~\ref{conj:5.2} and \ref{conj:3.3} are resolved completely.

Note that for a fixed color-critical graph $F$ with \( \chi(F) = r + 1 \geqslant 4 \), the graph in $\mathcal{G}(m, F)$ with maximum spectral radius is $r$-partite. Motivated by problem~\ref{pb:5.3}, it is naturally to consider the following problems. 
\begin{problem}\label{pb:5.4} 
Let \( F \) be a color-critical graph with \( \chi(F) = r + 1 \geqslant 4 \). For all sufficiently large $m$, what is the extremal graph with maximal spectral radius over all non-$r$-partite $F$-free graphs on $m$ edges?
\end{problem}
\begin{problem}
Solved Problem~\ref{pb:5.4} for complete graphs, generalized book graphs and even wheels.
\end{problem}
\section*{\normalsize Declaration of competing interest}
The authors declare that they have no known competing financial interests or personal relationships that could have influenced the work reported in this paper.

\section*{\normalsize Acknowledgement}
Minjie Zhang financially supported by the Open Research Fund of Key Laboratory of Nonlinear Analysis \& Applications (CCNU), Ministry of Education of China (Grant No. NAA2025ORG010).

\section*{\normalsize Data availability}
No data were utilized for the research that is described in this article.

\end{document}